\documentclass[preprint,12pt]{elsarticle}

\usepackage{graphicx}% Include figure files
\usepackage{dcolumn}% Align table columns on decimal point
\usepackage{bm}
\usepackage{epsfig}
\usepackage{amssymb}
\journal{Physica A}

\begin{document}

\begin{frontmatter}


%%\title{\tnoteref{label1}}

\title{Fixed points and stability in the two-network frustrated Kuramoto model}


\author[label1]{Alexander C. Kalloniatis}
\ead{alexander.kalloniatis@dsto.defence.gov.au}
\author[label1]{Mathew L. Zuparic}
\ead{mathew.zuparic@dsto.defence.gov.au}
%\address[label1]{Departamento de Fisica, Universidade Federal de Sao Carlos, 
%Sao Carlos, Brazil}
\address[label1]{Defence Science and Technology Group, Canberra,
ACT 2600, Australia}

\begin{abstract}
We examine a modification of the Kuramoto model for phase oscillators coupled on a network. Here, two populations of 
oscillators are considered, each with different network topologies, internal and 
cross-network couplings and frequencies. 
Additionally, frustration parameters for the interactions of the cross-network phases are introduced. This may 
be regarded as a model of competing populations: internal to any one network {\it phase} synchronisation is a 
target state, while externally one or both populations  seek to {\it frequency} synchronise to a phase in relation to 
the competitor. We conduct fixed point 
analyses for two regimes: one, where internal phase synchronisation occurs for each population with the potential for
instability in the phase of one population 
in relation to the other; the second where one part of a population remains fixed in phase in relation to the other population, but
where instability may occur within the first population leading to `fragmentation'. We compare analytic results to 
numerical solutions for the system at various critical thresholds.
\end{abstract}


\end{frontmatter}

% \linenumbers

%% main textzzzzzzzzzzzzzzzzzzzzzzzzzzzzzzzzzzzzzzzzzzzzzzzzzzzzzzzzzzzzzzzzzzzzzzzzzzzzzzzzzzzzzzzzzzzzzzzzzzzzzzzzzzzzzzzzzzzz

\section{Introduction}
Dynamical processes on networks remain an ongoing area of research across complex systems of vastly different
manifestations: physical, chemical, biological and social/organisational. 
The Kuramoto model of phase oscillators \cite{Kura1984} on a general network is the most paradigmatic of formulations allowing for
a myriad of variations depending on the specific area of application; for reviews see \cite{ABS2005,DGM2008,ADKMZ2008,DorBull2014}.
Among these, two variations of the Kuramoto model have recently become of interest in the literature: firstly, 
the {\it multi-network} formulation \cite{Bocc2014,MKB2004,BHOS2008,KNAKK2010,SkarRes2012}, 
where the entities being modelled on the network nodes may fall into quite different but nonetheless 
interacting populations (for example distinct species of
organisms) each with very different network characteristics; and secondly,
the {\it frustrated} Kuramoto model, where phase shifts are introduced in the interaction between adjacent oscillators such that the local interaction is not phase 
but frequency synchronising \cite{Cool2003,NVCDGL2013,KirkSev2015}. In this paper we consider an intersection of these two variations - a dual-network model, with frustrations introduced in the cross-network interaction.
Our interest in such a model stems from our adaptation of the concept of adversarial or competitive interaction in social systems
where two (or more) populations share a rivalry 
(for example, firms for market share \cite{Uzz97}, political or religious parties for membership \cite{AYW2011}, or 
popular opinion \cite{GACSM2014}). Note that the Kuramoto model has been adapted to social/human systems for several contexts:
rhythmic applause \cite{Ned2000}, opinion dynamics \cite{Plu2006} and human-robot musical performance \cite{Miz2010}.
There are also applications of coupled oscillator models on networks in which clustered populations are of interest, such as
in the work of \cite{Pec2014}.
Our application of the Kuramoto model is to the decision-making process because of its
essential {\it cyclicity}: for example the Perceptual Cycle model  of  Neisser \cite{Neis76} is based on a flow of an entity perceiving
external data, cognitively processing that data, making decisions about future actions, and undertaking them thereby depositing new
data into the external environment serving as an input into
other entities' process. In this respect then, the frustration in the cross-interaction may represent the aim of one group to be `a step ahead' of the
decision-making of competitors.
Inspired by this, the variation of the Kuramoto model we propose, for which we solve for fixed points, is close in spirit to the two network models of \cite{MKB2004,BHOS2008,KNAKK2010}
which also include frustrations - however we consider finite general networks as opposed to a coupling of two complete graphs. 

For a general unweighted 
undirected network described by an adjacency matrix $\cal A$, the Kuramoto model
is expressed by the coupled differential equations:
\begin{eqnarray}
{\dot \theta}_i = \omega_i + \sigma \sum_{j=1}^N {\cal A}_{ij} \sin(\theta_j-\theta_i), \;\; 1 \le i \le N
\end{eqnarray}
where $\theta_i$ is the phase angle for node $i$, $\omega_i$ is an intrinsic frequency associated 
with the node, and $\sigma$ a coupling constant.
The basic behaviour of this system has been well explored: for sufficiently high coupling, all phases 
may eventually synchronise and approximately phase lock: $\theta_i\approx \theta_j \ \forall i,j$
(there are exceptions though, see \cite{Tay2012}).

For two-network systems the equations are doubled and an additional cross network interaction is
introduced. Each network represents a distinct population, collection or organisation of entities described by the phase at the node.
To distinguish between the two competing populations we shall refer to them (based on competitive decision making approaches) respectively as 'Blue' and 'Red',
and the entities at nodes of the respective networks as agents. 
Analysing such a model for fixed points, 
Lyapunov stability/instability and chaos allows us to pose questions such as: 
How should agents be linked to each other? How quickly should linked agents respond to each other 
compared to responsiveness to changes in their specific competitor? How much diversity in frequency can be afforded by any one
population in light of its competitive interactions? And how can the aim for internal phase synchronisation be balanced against that
for frequency synchronisation (being 'ahead') of the competitor. We show how the model can be unpacked to answer such questions.

\begin{figure}[t]
\begin{center}
\includegraphics[height=5.7cm]{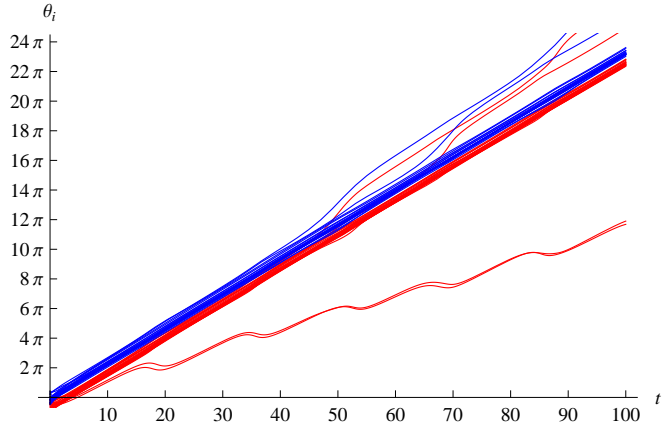}
\caption{
Typical example of behaviour of  phases $\theta_i$ as functions of time for the model introduced in this paper in which there are two populations
of phases `blue' and `red', where blue entities interact at some coupling strength and on some
network while red entities may have a different coupling strength and a different network, and some blue and red entities interact with each other
at yet another coupling strength and network. The splitting of the red lines, with two dropping periodically in relation to the main group of lines, is indicative of fragmentation of the
red population.}
\label{phases}
\end{center}
\end{figure}

While there are analogues of our approach in network control theory with a master network controlling 
another network, also known as `outer synchronisation'
\cite{AshMig2013},
in our case the two networks stand in a symmetric relationship to each other - each may be seen as seeking to control the other. 
A plot of individual $\theta_i$ as functions of time illustrates our focus. In Fig.\ref{phases} we show blue and red $\theta_i$
which represent our two types of entities. Different networks connect blue and red internally, with a whole other network for blue to red;
different couplings apply for entities of the same group or between the groups. The details behind this example are not important.
But we observe that most blue entities are locked slightly ahead of the majority of the red population, while two blue entities are detaching from and then reconnecting with their groups. When lines diverge they reach $2\pi$ difference from each other, namely they rejoin on the other side of the unit circle,
and another two red entities completely detach from their group, rejoin, and so on.

It is important
to note that our intent is {\it not} to seek prediction of how such network connected entities behave in general for low couplings. Such a regime would undermine control.
We seek, rather, to determine thresholds before `control' is lost, so that one population may achieve internal phase
synchronisation and external frequency synchronisation with respect to their competitor - such as the main blue population in Fig.\ref{phases}.
Within this we seek conditions such that the threshold for the type of fragmenation seen in the low-running red entities in  Fig.\ref{phases}
may be avoided. We derive such thresholds from linearisation and classical fixed point conditions. 
Such thresholds may be seen as bounding a `strategy' (for example, for competitive decision making) for one population in relation to the other.

In the next section we detail the model and then analyse its 
behaviour close to fixed points for the two populations locking internally and then for one population fragmenting into two parts. 
We then illustrate our results with numerical examples.
The paper concludes with a summary and future directions of the work.

\section{The frustrated two-network model: 'Blue vs Red'}
Consider $N$ Blue agents in an internally connected network given by an adjacency matrix ${\cal B}_{ij},\ ( i,j=1,\dots,N \in {\cal B} )$, and $M$ Red agents connected
in a network given by the adjacency matrix ${\cal R}_{ij},\ (i,j=N+1,\dots, N+M \in {\cal R})$,
namely matrices with entries one if $i$ and $j$ are connected by a link and zero otherwise.
Associated with each Blue agent $i \in {\cal B}$, is a phase $\beta_i$ giving the position in a limit cycle; similarly $\rho_i$ is the position in the limit cycle of a Red agent
$i\in{\cal R}$.
The 'Blue vs Red model' is given by the system of equations:
\begin{eqnarray}
\dot{\beta}_i = \omega_i + \sigma_B \sum_{j\in {\cal B}}{\cal B}_{ij} \sin(\beta_j-\beta_i) 
+ \zeta_{BR} \sum_{j \in {\cal R}} {\cal M}_{ij} \sin(\rho_j+\phi-\beta_i) ,\;\;  i \in {\cal B} \label{Blue-eq} \\
\dot{\rho}_i = \nu_i + \sigma_R \sum_{j \in {\cal R}} {\cal R}_{ij} \sin(\rho_j-\rho_i) 
+ \zeta_{RB} \sum_{j \in {\cal B}} {\cal M}_{ij} \sin(\beta_j+\psi-\rho_i), \;\;  i\in {\cal R}.  \label{Red-eq}
\end{eqnarray}
The matrix ${\cal M}_{ij}$ is block off-diagonal of dimension $N+M$ representing the network of links of Blue to Red 
and {\it vice versa}. The non-zero blocks of ${\cal M}$ will be denoted by the $N\times M$ and $M\times N$ matrices $A^{(BR)}$ and $A^{(RB)}$ respectively:
\begin{eqnarray}
{\cal M} = \left(\begin{array}{cc}
0 & A^{(BR)} \\
A^{(RB)} & 0 \end{array} \right).  
\end{eqnarray}
The quantities $\omega_i, \nu_i$ 
represent frequencies of individual Blue and Red agents respectively; these may be fixed or drawn from some probability distribution.
The frustration parameters $\phi,\psi$ represent the degree to which Blue and Red, respectively, seek to be ahead of their competitor's limit cycle. Indeed, by
representing this interaction as a shifted sine function we are ensuring that Blue, respectively Red, seek to stay $\phi$, respectively $\psi$, ahead of the other
but do not continue to accelerate; in other words there is no advantage to continually lap the adversary, otherwise agents are engaging in 
increasingly disconnected limit cycles. 
Finally, $\sigma_B, \sigma_R, \zeta_{BR}, \zeta_{RB}$ are coupling constants, respectively, for intra-Blue, intra-Red, Blue to Red and Red to Blue. 
Asymmetry between Blue and Red potentially exists both in the coupling constant and the network: $\zeta_{BR}$ need not equal $\zeta_{RB}$ and a network link 
from Blue agent $i$ to Red agent $j$ need not imply a link from Red agent $j$ to Blue agent $i$, so that
the cross-network matrix ${\cal M}_{ij}$ need not be symmetric. 
Note, however, the overall symmetry between Blue and Red in Eqs.(\ref{Blue-eq},\ref{Red-eq}): neither side can be regarded as a master controller, though each seeks to 
assert control over the other if the nonlinear dynamics permits. In that respect, as we said at the outset, it is more appropriate to view this
not in terms of a controller but of a game between two teams of, respectively, $N$ and $M$ players. In this sense, the choice of
frustration parameters $\phi,\psi$ may be referred to as a `strategy':
success for agents in these systems means the simultaneous achievement of two activities: synchronising internal 
limit cycles and staying as close as possible to a fixed phase ahead of the limit cycle of the competitor.

\section{Fixed Point Analysis}
\subsection{Defining the fixed point for internal locking}
The system Eqs.(\ref{Blue-eq},\ref{Red-eq}) in general can only be solved numerically. However, since both sides 
may deem internal phase synchronisation of limit cycles ideal we may explore the regime of a fixed point given by the following conditions:
\begin{eqnarray}
\beta_i  =  B + b_i, \;\; \rho_i  =  P + p_i. \label{fixpoint1}
\end{eqnarray}
The variables $B, b_i, P, p_i$ are all time dependent, but $b_i,p_i$ are `small' fluctuations, and hence we assume $b_i^2 \approx 0$,  $p^2_j \approx 0$, $\forall \;\{i,j\} \in \{{\cal B},{\cal R}\}$.
The variables $B,P$ thus represent the centroid of the Blue and Red phases respectively
\begin{eqnarray}
B= {1\over{N}} \sum_{i\in {\cal B}} \beta_i , \;\; P= {1\over{M}} \sum_{i \in {\cal R}} \rho_i \label{BPdef}.
\end{eqnarray}
The conditions in Eq.(\ref{fixpoint1}) follow analysis of the pure Kuramoto model in \cite{Kall10}, using the idea of the Master Stability Function for network coupled
dynamical systems by \cite{Pec-Car98}.
The difference between the sides' global phases is:
\begin{eqnarray}
\alpha \equiv B-P  \label{fixpoint-alph}
\end{eqnarray}
We shall soon examine the conditions under which $\alpha$ may be constant, but initially - even with Blue and Red internally phase locking -
there may be time dependence within which one or more members of a population loops around the other on the unit circle. 
Thus, as a consequence of any one term in the sums in Eq.(\ref{BPdef}) incrementing by $2\pi$, $\alpha$ should be
treated as modulo $2\pi/N$ (from $B$) or $2\pi/M$ (from $P$). However, because of the linearisation no fluctuation $b_i, p_i$ can develop to represent
one or few members of a population breaking off from the whole - we later allow for a separate cluster in the {\it ansatz} - so that 
 $-\pi \leq \alpha \leq \pi$ is the natural range within this {\it ansatz}. We refer to the phase locking within Blue, $\beta_i \approx \beta_j\; \forall\; i,j \in {\cal B}$, 
and Red $\rho_i \approx \rho_j\; \forall \; i,j \in {\cal R}$, respectively, as {\it local locking}, and the
phase locking of Blue with respect to Red, $\beta_i \approx \rho_j \; \forall \;\{i,j\} \in \{{\cal B},{\cal R}\}$, as {\it global locking}.
Later we shall examine the case of one population fragmenting, but this serves to establish our notation and methods.

\subsection{Linearisation for internal locking}
We now linearise the system around conditions Eqs.(\ref{fixpoint1},\ref{fixpoint-alph}) using the smallness of 
$b_i,p_i$ at this point, $\phi,\psi,\alpha$ being arbitrary. We use the following 
Taylor approximations of the interaction terms:
\begin{eqnarray}
\sin(\beta_j-\beta_i) \approx b_j-b_i, \;\; \sin(\rho_j-\rho_i) \approx p_j-p_i \nonumber \\
\sin(\rho_j+\phi-\beta_i) \approx \sin(\phi-\alpha) + (p_j-b_i) \cos(\phi-\alpha) \nonumber \\
\sin(\beta_j+\psi-\rho_i) \approx \sin(\psi+\alpha) + (b_j-p_i) \cos(\psi+\alpha). \label{approx}
\end{eqnarray}
Inserting Eqs.(\ref{approx}) into Eqs.(\ref{Blue-eq},\ref{Red-eq}), the resulting system can be cast in the form
\begin{eqnarray}
\dot{\vec{v}} + \dot{\vec V} = \vec{\Omega} - {\mathcal L} \vec{v}, \label{linsys}
\end{eqnarray}
where
\begin{eqnarray}
\vec{v} = \left( \begin{array}{c} {\hat b} \\ {\hat p} \end{array} \right),\;\; \vec {V}  =\left( \begin{array}{c} B {\hat 1}_B \\ P {\hat 1}_R \end{array} \right), \;\;
\vec{\Omega} = \left( \begin{array}{c} \hat\omega+ \zeta_{BR}\sin(\phi-\alpha) {\hat d}^{(BR)} \\ \hat\nu+ \zeta_{RB}\sin(\psi+\alpha) {\hat d}^{(RB)} \end{array}
\right) \label{fvec} 
\end{eqnarray}
and
\begin{eqnarray}
&&{\mathcal L} = \label{genLapl}\\
&&\left( \begin{array}{cc} \sigma_B L^{(B)}+ \zeta_{BR}\cos(\phi-\alpha) D^{(BR)} & - \zeta_{BR}\cos(\phi-\alpha) A^{(BR)} \\
- \zeta_{RB} \cos(\psi + \alpha) A^{(RB)} & \sigma_R L^{(R)} + \zeta_{RB} \cos(\psi + \alpha) D^{(RB)} \end{array} \right). \nonumber
\end{eqnarray}
We have adopted the vector notation ${\vec v}$ for $N+M$ dimensional objects across the full Blue-Red system, while reserving 
`hat' quantities such as ${\hat e}$ for $N$ or $M$ dimensional 
vectors relevant to the Blue or Red networks. 
The vectors ${\hat 1}_B, {\hat 1}_R$ have $N$, respectively $M$, components all value one. The quantities ${\hat d}, D$ and $L$ are explained in the following.

The matrices $L$ inside the blocks of $\cal L$ constitute a {\it graph Laplacian}. To explain this we introduce some basic notions
of graph theory from \cite{Boll98} using the example of the Blue population. The degree of each Blue agent $i$ is the number of links from $i$ to other Blue agents,
\begin{eqnarray}
d^{(B)}_i = \sum_{j\in {\cal B}} {\cal B}_{ij}. \nonumber
\end{eqnarray}
As an $N$ dimensional vector this is written as ${\hat d}^{(B)}$. We then introduce a matrix $D^{(B)}$ whose only non-zero elements are the degrees $d^{(B)}_i$ running along the diagonal
\begin{eqnarray}
{\cal D}^{(B)}_{ij} = d^{(B)}_i \delta_{ij}. \nonumber 
\end{eqnarray}
The Laplacian for the Blue population is then
\begin{eqnarray}
L^{(B)}_{ij} = {\cal D}^{(B)}_{ij} - {\cal B}_{ij}. \nonumber 
\end{eqnarray}
Equivalent definitions apply to the Red population leading to the graph Laplacian for Red, $L^{(R)}$.

The matrices $D^{(BR)}, D^{(RB)}$ are similar to the degree matrix but encode the Blue to Red links. 
For example, we define a degree with which a Blue node $i$ connects to Red agents,
\begin{eqnarray}
d^{(BR)}_i = \sum_{j\in {\cal R}} {\cal M}_{ij}, \nonumber
\end{eqnarray}
which is written as the vector ${\hat d}^{(BR)}$ in Eq.(\ref{fvec}).
The degree matrix can then correspondingly be formed:
\begin{eqnarray}
D^{(BR)}_{ij} = d^{(BR)}_i \delta_{ij},\;\; i,j \in {\cal B}  \nonumber 
\end{eqnarray}
The same considerations enable definition of $D^{(RB)}$. 

\subsection{The free system}\label{FREESYS}
We would like to decouple the linear system of equations in Eq.(\ref{linsys}) by expanding in some suitable set of eigenvectors.
Here the spectral properties of Laplacians become a powerful tool, see \cite{Boll98}. 
However, as ${\cal L}$ is not generally symmetric, its left and right eigenvectors are not
the hermitian conjugates of each other. The asymmetry in ${\cal L}$ is generated by the cross-network interactions. 
For $\zeta_{BR}=\zeta_{RB}=0$ these problems disappear:
\begin{eqnarray}
{\cal L}={\cal L}_0 \equiv \left( \begin{array}{cc} \sigma_B L^{(B)} & 0 \\
0 & \sigma_R L^{(R)} \end{array}\right).  \label{freeLapl}
\end{eqnarray}
To orient ourselves through the problem and to set up notation for the interacting situation
it is useful to analyse this case carefully. Using our knowledge of the spectral properties of graph Laplacians, a
spanning set of orthonormal eigenvectors for the super-Laplacian in the $N+M$ dimensional space is:
\begin{eqnarray*}
{\vec e}^{(r)}= \left\{ \begin{array}{lll} \left(  \begin{array}{c} {\hat e}^{(B,r)}  \\ 0 \end{array}\right) &\textrm{for}&    r=0, \dots, N-1 \\
  \left( \begin{array}{c} 0 \\ {\hat e}^{(R,r)}  \end{array}\right) &\textrm{for}&   r=N, \dots, N+M-1  \end{array} \right.
\end{eqnarray*}
Observe how, consistent with starting our index from zero for the first (Blue population Laplacian) zero mode, we have reserved the index value
$r=N$ for the Red population Laplacian zero mode. 
These (orthonormal) zero-modes are given through the identification ${\hat e}^{(B,0)}=\hat{1}_B/\sqrt{N}, {\hat e}^{(R,N)}=\hat{1}_R/\sqrt{M}$.
We have seen that the vector ${\vec V}$ is given in terms of these zero eigenvectors, thus we expand the fluctuations ${\vec v}$ in the `normal' modes:
\begin{eqnarray}
{\vec v} = \sum_{r \neq 0, N}^{N+M-1} x_r {\vec e}^{(r)}. \label{normalex}
\end{eqnarray}
Projecting Eq.(\ref{linsys}) from the left on the two zero modes, respectively normal modes, we obtain:
\begin{eqnarray}
{\dot B} = {\bar \omega}; \ r=0, \;\; {\dot P} = {\bar \nu}; \ r=N, \;\;{\dot x}_r = f_r - \lambda_r^{(0)} x_r; \ r\neq 0,N   \label{freesys}
\end{eqnarray}
where $f_r= {\vec \Omega} \cdot {\vec e}_r $ is the projection of ${\vec \Omega}$ onto the $r$-th eigenvector, $\bar\omega, \bar\nu$ represent the mean
frequencies within the blue, respectively red, populations and 
\begin{eqnarray}
\lambda_r^{(0)}= \left\{ \begin{array}{lll} \sigma_B \lambda_r^{(B)}&\textrm{for}& r=1,\dots,N-1  \\
 \sigma_R \lambda_r^{(R)}&\textrm{for}&  r=N+1,\dots,N+M-1 \end{array}\right. \label{e-val-0}
\end{eqnarray}
are the eigenvalues. We may draw two conclusions out of this. 

Firstly, the condition for the free system to have Blue agents lock a fixed phase angle in relation to Red agents amounts to 
${\dot \alpha}={\dot B}-{\dot P}=0$, thus
\begin{eqnarray}
\bar \omega = \bar \nu.\label{freelock}
\end{eqnarray}
This is intuitively sound since our analysis presupposes that Blue and Red networks have locally locked, and there is no interaction between them
so the appearance of global locking is fortuitous.
Note that even if both Blue and Red frequencies are drawn from the same statistical distribution, for finite $N,M$ this condition will not generally be
satisfied, in other words $\bar \omega \neq \bar \nu$, so that one system will run ahead of the other, according to which has greater average frequency,
with time dependent $\alpha$ global phase difference. In fact, the equation for $\dot\alpha$ can be obtained directly
by projection of Eq.(\ref{linsys}) with 
\begin{eqnarray}
{\vec {\Pi}} \equiv 
\left( \begin{array}{c}
{{\hat{1}_B}\over{N}} \\
-{{\hat{1}_R}\over{M}} \end{array}\right).
\end{eqnarray}
This is evidently a linear combination of zero modes of ${\cal L}_0$ so that fluctuations ${\vec v}$ decouple from
the dynamics of $\alpha$. We will see that this no longer applies for the interacting system.

A second conclusion is that the system is absolutely stable due to the positivity of the normal mode spectrum - a consequence of the property of the separate Blue and Red Laplacians -
manifesting the absence of cross-network interactions to challenge internal synchronisation. Also, since this analysis is predicated on small fluctuations, 
$x_r(t)$ must remain small, giving, from Eq.(\ref{freesys}),
\begin{eqnarray}
\left|f_r/\lambda_r^{(0)} \right| \ll 1. \label{critcond}
\end{eqnarray}
With $\lambda_r^{(0)}=\sigma_B \lambda_r^{(B)},\sigma_R \lambda_r^{(R)}$ (Eq.(\ref{e-val-0})) this implies a lower bound on the coupling in terms of the frequency range
residing in $f_r$. As we discuss in \cite{Zup2013}, this is only a weak constraint there being no discrete scale which says how small is sufficient,
and its violation does not generate an instability. However, it does provide a guide to reasonable values of coupling to permit 
internal synchronisation of agents within the respective populations. 

\subsection{Interacting system - zero modes}
For the interacting system left eigenvectors cease to be the conjugate of right eigenvectors. This means that the projector for the equation for $\dot{\alpha}$ ceases to 
annihilate against the super-Laplacian $\cal L$ and $\dot{\vec v}$ remains in the equation.
Nevertheless, we persist with the free system basis - accepting that it will not diagonalise ${\cal L}$. In this case the dynamical equation for the angle between the centroids of the two populations is
\begin{eqnarray}
\dot{\alpha} = {\vec \Pi} \cdot \left(  {\vec\Omega} - {\cal L} {\vec v}  \right). \label{totalalpha}
\end{eqnarray}
We see now that the approximate equation $\dot{\alpha} \approx \Pi \cdot  {\vec\Omega}$ only holds if we suppress
the normal mode fluctuations. In this approximation the solution is:
\begin{eqnarray}
\dot{\alpha} \approx  \bar{\omega} - \bar{\nu}  + S(\phi,\psi) \cos\alpha - C(\phi,\psi) \sin\alpha, \label{alphadot-approx}
\end{eqnarray}
where
\begin{eqnarray*}
&C(\phi,\psi)  \equiv  \frac{d^{(BR)}_T \zeta_{BR} \cos\phi}{N} + \frac{d^{(RB)}_T \zeta_{RB} \cos\psi}{M},\;\; S(\phi,\psi)  \equiv  \frac{d^{(BR)}_T \zeta_{BR}\sin\phi}{N}  -  \frac{ d^{(RB)}_T \zeta_{RB} \sin\psi   }{M}
\end{eqnarray*}
and $d^{(BR)}_T\equiv \sum_{i,j}A^{(BR)}_{ij}$ is the total degree between Blue and Red.
Eq.(\ref{alphadot-approx}) is straightforwardly integrated, with the help of the Weierstrass substitution 
$\alpha = 2 \arctan(\alpha')$ and standard trigonometric formulae, giving the time-dependent solution:
\begin{eqnarray}
\alpha(t) =
2 \arctan \left\{ \frac{ C(\phi,\psi) - \sqrt{{\cal K}(\phi,\psi)}
\tanh\left(\frac{t + {\rm{const}}}{2} \sqrt{{\cal K}(\phi,\psi)} \right)}{\bar{\omega} - \bar{\nu} -  S(\phi,\psi)} \right\}
\label{t-alpha}
\end{eqnarray}
where `const' is the integration constant given through initial conditions, and
\begin{eqnarray}
{\cal K}(\phi,\psi) \equiv C(\phi,\psi)^2 + S(\phi,\psi)^2- (\bar{\omega} - \bar{\nu})^2. \label{specialK}
\end{eqnarray}

If couplings are such that $\alpha$ reaches a steady-state value then 
$\dot{\alpha}=0$ and $\alpha$ is more easily determined from the right hand side of Eq.(\ref{alphadot-approx}):
\begin{eqnarray}
 \bar{\omega} - \bar{\nu}  + S(\phi,\psi) \cos\alpha - C(\phi,\psi) \sin\alpha \approx 0. 
\label{alpha-const-sol}
\end{eqnarray}
We observe corrections in Eq.(\ref{alpha-const-sol}) compared to the free-system result Eq.(\ref{freelock}) coming from interactions.
This equation is readily solved for $\sin\alpha$:
\begin{eqnarray}
\sin\alpha &=& \frac{(\bar{\omega} - \bar{\nu}) C(\phi,\psi)  \pm S(\phi,\psi) \sqrt{
{\cal K(\phi,\psi)}}}{C(\phi,\psi)^2+S(\phi,\psi)^2}
. \label{alphasol}
\end{eqnarray}

We highlight a number of properties in Eq.(\ref{alphasol}). Firstly, there will be complex solutions
when ${\cal K}<0$. This coincides with the point in Eq.(\ref{t-alpha}) where the time
dependence goes from hyperbolic to trigonometric tangent (namely, from a plateau to oscillatory behaviour).
Secondly there may be no solutions when the right hand side
of Eq.(\ref{alphasol}) exceeds $\pm 1$. 
Both are indicators that the solution Eq.(\ref{t-alpha}) is varying for all $t$, and therefore that there does not exist a fixed point.

This boundedness of constant solutions
$\alpha$, and the dependence of the solution on trigonometric functions of
the frustration parameters $\phi$ and $\psi$ means that the preferred fixed point
of the dynamical system is nonlinear in the parameters defining a population's strategy.
Stated simply, there may be {\it diminishing returns in Blue seeking to be too far ahead of Red's limit cycle}.
We can state this as an optimisation problem. Take the solution to Eq.(\ref{alphasol}) as the function $\alpha(\phi,\psi)$, of
the frustration parameters. This may be regarded as an {\it objective function} determining the optimal strategy for Blue or Red.
{\it Blue seeks a strategy $\phi$, given structures, couplings and Red's strategy $\psi$, that maximises $\alpha(\phi,\psi)$.
Conversely, Red seeks a strategy $\psi$, given structure, couplings and Blue's strategy $\phi$, that
maximises $-\alpha(\phi,\psi)$.} We show the result of such considerations for a specific numerical example later in
the paper.

\subsection{Conditions on the Spectrum}
For the interacting system now we write the eigenvalue-eigenvector problem for the Blue population as
\begin{eqnarray*}
\sum_{j \in {\cal B}} L^{(B)}_{ij} e^{(B,r)}_j &=& {\lambda}^{(B)}_r e^{(B,r)}_i  
\end{eqnarray*}
where $r$ labels the $N$ independent eigenvectors ${\hat e}$ and discrete eigenvalues $\lambda$. 
Laplacians so structured enjoy a positive semi-definite spectrum, with the degeneracy of zero eigenvalues, which we label ${\lambda}_0=0$,
corresponding to the number of disconnected components the network breaks up into without removing any of the existing links or nodes.
The corresponding un-normalised eigenvector is that with only the values one in its entries: $e^{(B,0)}_i \propto 1 \  \forall \  i=1,\dots,N$. We
have already encountered these vectors in the definition of ${\vec V}$. Thus, $e^{(B,0)}=\frac{1}{\sqrt{N}}{\hat 1}_B$, and
$e^{(R,0)}=\frac{1}{\sqrt{M}}{\hat 1}_R$.
In the following, we consider Blue and Red networks that each have only one component so that there is only one zero mode for each of $L^{(B)},L^{(R)}$.
The lowest non-zero eigenvalue, for example for Blue, ${\lambda}^{(B)}_1>0$ is known as the algebraic connectivity, 
\cite{Fied73}, encoding properties of the extent to which nodes are strongly or weakly linked.
For example, networks with $\lambda^{(B)}_1<1$ are easily disconnected into separate components with the removal of a small number of Blue nodes or links.
The components of the corresponding eigenvector $e^{(B,1)}_i$ behave like a step function in $i$, with jumps at node values coinciding with the
most poorly connected nodes in the network (those whose removal decouples the network), 
as illustrated in the work by \cite{Ding-et-al01}.
The same statements apply for the Red population.

As mentioned, ${\cal L}$ is not symmetric - except when $\zeta_{BR} \cos(\phi-\alpha)=\zeta_{RB} \cos(\psi+\alpha)$ and every Blue node that links to a Red node in turn is linked to by the same Red node. Therefore, even though its matrix elements
are real, ${\cal L}$ may have complex eigenvalues. Unfortunately, the Gershgorin theorem \cite{Gersh31} only provides bounds on discs within which eigenvalues must lie. 
An alternative uses a trick employed in the proof
of a lemma by Taylor \cite{Tay2012} which identifies positive-definiteness of the super-Laplacian ${\cal L}$:
for every $\vec{z}$ that $\vec{z}^T {\cal L} \vec{z}\geq 0$. Choosing the free-system eigenvectors
$\vec{e}^{(0)}, \vec{e}^{(N)}$, which trivially annihilate under ${\cal L}_0$ and against the off-diagonal parts of ${\cal L}$, we obtain 
from the remaining diagonal elements the two conditions
\begin{eqnarray}
\cos(\phi-\alpha) \geq 0, \;\;\cos(\psi+\alpha) \geq 0.
\label{stab-2clust}
\end{eqnarray}
Note that for non-symmetric matrices the {\it absence} of positive semi-definiteness is necessary but not sufficient for there to
be negative eigenvalues (for example the matrix $((1,-5),(0,1))$ is not positive definite as seen using the vector $(1,1)$, but has 
eigenvalue $+1$). The same inequality can be obtained from the Gershgorin theorem. However, when $\zeta_{BR}\cos(\phi-\alpha)=\zeta_{RB}\cos(\psi+\alpha)$,
where the super-Laplacian is symmetric, then the above trick gives necessary {\it and sufficient} conditions for positive eigenvalues
 - whereas the Gershgorin theorem still only gives a left most point to an interval which {\it contains} an eigenvalue.

\subsection{Red fragmentation}
\label{3cent}
We now consider the case that Red may possibly fragment, with those parts of it interacting with Blue agents remaining locked with them.
The case of Blue fragmenting is entirely symmetric and may be extracted from the following analysis by swapping labels.
We partition Red into ${\cal R}_1$ and ${\cal R}_2$ with $M_1$ and $M_2$ nodes each, $M_1+M_2=M$.
We will later choose ${\cal R}_1$ as those Red agents interacting with Blue agents, but for completeness the full system of equations
after this partition reads:
\begin{eqnarray*}
\dot{\beta}_i = \omega_i - \sigma_B \sum_{j \in {\cal B}} {\cal B}_{ij} \sin \left(\beta_i-\beta_j \right) 
- \zeta_{BR} \sum_{j\in {\cal R}_1} {\cal M}_{ij} \sin(\beta_i-\phi-\rho_j) \\
- \zeta_{BR} \sum_{j\in {\cal R}_2} {\cal M}_{ij} \sin(\beta_i-\phi-\rho_j) ,\;\;  i \in {\cal B} \\
\dot{\rho}_i = \nu^{(a)}_i - \sigma_R \sum_{j \in {\cal R}_1} {\cal R}_{ij} \sin \left(\rho_i-\rho_j \right)- \sigma_R \sum_{j \in {\cal R}_2} {\cal R}_{ij} \sin \left(\rho_i-\rho_j \right) \\
 - \zeta_{RB} \sum_{j\in {\cal B}} {\cal M}_{ij} \sin(\rho_i-\psi-\beta_j) ,\;\;  i \in {\cal R}_a ,& a=\{1,2\} 
\end{eqnarray*}
The index $a$ labels the parts of ${\cal R}$. 
Analogously to our analysis of internal locking, we now impose the fixed point conditions
\begin{eqnarray}
\beta_i = B + b_i,\; i \in {\cal B},\;\; \rho_j = P_l + p^{(a)}_j, \;  j \in {\cal R}_a, \label{fixpoint2} 
\end{eqnarray}
and define the three `centroid' phases
\begin{eqnarray}
 B-P_1\equiv \alpha_{B R_1},\;\;B-P_2 \equiv \alpha_{B R_2} \;\; P_1 - P_2 \equiv \alpha_{R_1 R_2}. \label{3-cent-def}
\end{eqnarray}
These may be written in terms of the $\beta_i$ and $\rho_i$ using analogues of Eq.(\ref{BPdef}).

Applying these expansions and definitions, the linearised equations Eq.(\ref{linsys}) become $3\times 3$ dimensional, where
the vectors and matrices there are replaced by `primed' quantities:
\begin{eqnarray}
\vec{v'} =\left( \begin{array}{c} {\hat b} \\ {\hat p^{(1)}} \\  {\hat p^{(2)}} \end{array} \right),  \ \
\vec {V'} =\left( \begin{array}{c} B {\hat 1}_B \\ P_1 {\hat 1}_{R_1} \\  P_2 {\hat 1}_{R_2} \end{array} \right), \label{v'V'-defs} 
\end{eqnarray}
\begin{eqnarray}
\vec{\Omega}' = \left( \begin{array}{c} 
\hat\omega+ \zeta_{BR}\left\{ \sin(\phi-\alpha_{BR_1}) {\hat d}^{(BR_1)} + \sin(\phi-\alpha_{BR_2}) {\hat d}^{(BR_2)}\right\}
\\ 
\hat{\nu}^{(1)}+ \zeta_{RB}\sin(\psi+\alpha_{BR_1}) {\hat d}^{(R_1B)}-\sigma_R \sin (\alpha_{R_1 R_2})   {\hat d}^{(R_1 R_2)}
\\
\hat{\nu}^{(2)}+ \zeta_{RB}\sin(\psi+\alpha_{BR_2}) {\hat d}^{(R_2B)} + \sigma_{R}\sin(\alpha_{R_1 R_2}) {\hat d}^{(R_2 R_1)} 
\end{array}
\right)
\label{f'vec} 
\end{eqnarray}
and the super-Laplacian
\begin{eqnarray}
&&\!\!\!\!\!\!
{\mathcal L}' =\label{genLapl'}\\
&&\!\!\!\!\!\!\left( \begin{array}{ccc} 
\sigma_B L^{(B)}+{\cal V}_1
& - \zeta_{BR}\cos(\phi-\alpha_{BR_1}) A^{(BR_1)} 
&  -\zeta_{BR}\cos(\phi-\alpha_{BR_2}) A^{(BR_2)} 
\\
- \zeta_{RB} \cos(\psi+\alpha_{BR_1}) A^{(R_1 B)} 
&  \sigma_R L^{(R_1)} + {\cal V}_2 
&  -\sigma_R \cos(\alpha_{R_1R_2}) A^{(R_1 R_2)}
\\
- \zeta_{RB}\cos(\psi+\alpha_{BR_2}) A^{(R_2 B)} 
&
-  \sigma_R \cos(\alpha_{R_1R_2}) A^{(R_2 R_1)}
& \sigma_R L^{(R_2)} +  {\cal V}_3
 \end{array} \right)   \nonumber
\end{eqnarray}
We recognise in the diagonal elements the Laplacians for the disconnected sub-graphs ${\cal B}, {\cal R}_1, {\cal R}_2$, and
additional interaction-dependent contributions:
\begin{eqnarray*}
{\cal V}_1 =  \zeta_{BR}\cos(\phi-\alpha_{BR_1}) D^{(BR_1)} + \zeta_{BR}\cos(\phi-\alpha_{BR_2}) D^{(BR_2)} \\
{\cal V}_2 = \zeta_{RB} \cos(\psi + \alpha_{BR_1}) D^{(R_1B)} + \sigma_R \cos(\alpha_{R_1 R_2}) D^{(R_1 R_2)} \\
{\cal V}_3 = \zeta_{RB} \cos(\psi + \alpha_{BR_2}) D^{(R_2B)} + \sigma_R \cos(\alpha_{R_1 R_2}) D^{(R_2 R_1)} .
\end{eqnarray*}
Evidently now the free part of the super-Laplacian has three zero eigenvectors which are composed of the $N, M_1$ and $M_2$
dimensional vectors $\hat{1}_B, \hat{1}_{R_1}, \hat{1}_{R_2}$ all with entries equal to one.

In the presence of stability (conditions to be determined below) with small fluctuations, the dynamical equations for the centroids $\dot{\vec V}'={\vec\Omega}'$ lead to:
\begin{eqnarray}
\dot{B} = \bar{\omega} - \frac{\zeta_{BR} d^{(B R_1)}_T}{N} \sin \left( \alpha_{B R_1}-\phi \right) - \frac{\zeta_{BR}d^{(B R_2)}_T}{N} \sin (\alpha_{B R_2}-\phi)  \nonumber \\
\dot{P}_1 = \bar{\nu}^{(1)} -\frac{\sigma_{R} d^{(R_1 R_2)}_T}{M_1} \sin \left( \alpha_{R_1 R_2} \right)+ \frac{\zeta_{RB}d^{(R_1 B)}_T}{N} \sin (\alpha_{B R_1}+\psi) \nonumber \\
 \dot{P}_2 = \bar{\nu}^{(2)} + \frac{\sigma_{R} d^{(R_1 R_2)}_T }{M_2} \sin \left(\alpha_{R_1 R_2}\right)  + \frac{\zeta_{RB}d^{(R_2 B)}_T}{M_2} \sin (\alpha_{B R_2}+\psi) \label{3-cent-exp}
\end{eqnarray}
From the definitions of the three centroids given by Eq.(\ref{3-cent-def}) we can see that,
\begin{eqnarray}
\alpha_{B R_2}= \alpha_{B R_1}+\alpha_{R_1 R_2}.
\label{13=12+23}
\end{eqnarray}
Since the angles between the centroids are linearly dependent through Eq.(\ref{13=12+23}), one of the dynamical equations for the angles derived from Eq.(\ref{3-cent-def}) may be eliminated. Choosing $\alpha_{B R_2}$, we obtain the equations for the remaining angles between the centroids:
\begin{eqnarray*}
\dot{\alpha}_{B R_1} &=& \bar{\omega}- \bar{\nu}^{(1)} 
-\frac{\zeta_{BR} d^{(B R_1)}_T}{N} \sin \left( \alpha_{B R_1} -\phi \right)
- \frac{\zeta_{RB}d^{(R_1 B)}_T}{M_1} \sin (\alpha_{B R_1}+\psi)  \\ 
&&
+\frac{\sigma_R d^{(R_1 R_2)}_T}{M_1} \sin (\alpha_{ R_1 R_2})  
-\frac{\zeta_{BR}d^{(B R_2)}_T}{N} \sin (\alpha_{B R_1}+\alpha_{R_1 R_2}-\phi) \\
\dot{\alpha}_{R_1 R_2} &=& \bar{\nu}^{(1)}- \bar{\nu}^{(2)} - \frac{\zeta_{RB}d^{(R_2 B)}_T}{M_2} \sin (\alpha_{B R_1}+\alpha_{R_1 R_2}+\psi) \\
&&+\frac{\zeta_{RB} d^{(R_1 B)}_T}{M_1} \sin ( \alpha_{B R_1} + \psi ) -\sigma_{R}d^{(R_1 R_2)}_T\left( \frac{1}{M_1}+\frac{1}{M_2} \right) \sin (\alpha_{R_1 R_2}).  
\end{eqnarray*}

\subsection{Critical values and stability}
To be more specific now, it is natural to expect that with the interaction with Blue, the parts of Red most likely to fragment in relation to
each other are those Red agents that are interacting with Blue against those that are not. We choose then ${\cal R}_1$ to be
those interacting with Blue, so that $d^{(B R_2)}_T = d^{(R_2 B)}_T = 0$. 
Thus, at steady-state, $\dot{\alpha}_{B R_1}=\dot{\alpha}_{R_1 R_2}=0$, the system becomes, 
\begin{eqnarray*}
\chi_1 &=&C_1  \sin \left( \alpha_{B R_1} \right)- S_1  \cos \left( \alpha_{B R_1} \right) \\
\sin \left( \alpha_{R_1 R_2} \right) &=& \chi_2 + C_2  \sin \left( \alpha_{B R_1} \right) -S_2  \cos \left( \alpha_{B R_1} \right) ,
\end{eqnarray*}
where $\chi_1,\chi_2$ contain the frequencies, and $C_1,C_2,S_1$ and $S_2$ absorb all the terms multiplying cosines and sines
of $\phi$ and $\psi$:
\begin{eqnarray*}
&\chi_1 = \bar{\omega}-\bar{\nu}^{(2)}+\frac{M_1 }{M_2 }\left(  \bar{\omega}-\bar{\nu}^{(1)} \right),\;\; \chi_2 =
\frac{M_1 }{\sigma_R d^{(R_1 R_2)}_T}\left(  \bar{\nu}^{(1)} -\bar{\omega} \right) &\\
&C_1 = \frac{M \zeta_{BR} d^{(B R_1)}_T \cos\phi }{M_2 N}
+ \frac{\zeta_{RB} d^{(R_1 B)}_T\cos\psi}{M_2} , \;\; S_1 = \frac{M \zeta_{BR} d^{(B R_1)}_T\sin\phi }{M_2 N} 
- \frac{\zeta_{RB} d^{(R_1 B)}_T\sin\psi }{M_2}&  \\
&C_2= \frac{M_1\zeta_{BR}d^{(B R_1)}_T\cos \phi}{N \sigma_R d^{(R_1 R_2)}_T} +
\frac{\zeta_{RB}d^{(R_1 B)}_T\cos \psi}{ \sigma_R d^{(R_1 R_2)}_T}, \;\; S_2
= \frac{M_1\zeta_{BR}d^{(B R_1)}_T\sin \phi}{N \sigma_R d^{(R_1 R_2)}_T} -
\frac{\zeta_{RB}d^{(R_1 B)}_T\sin \psi}{ \sigma_R d^{(R_1 R_2)}_T}.&
\end{eqnarray*}
Solving for the sine of each of the angles we obtain,
\begin{eqnarray}
\sin(\alpha_{B R_1})&=& \frac{\chi_{1}C_{1} \pm S_{1} \sqrt{{\cal J}}}{C^2_{1} + S^2_{1}}\label{alphaBR_1} \\
\sin(\alpha_{R_1 R_2}) &=& \chi_2 +   \frac{\chi_{1}(C_{1}C_2+S_1 S_2)}{C^2_{1} + S^2_{1}}\pm \frac{(S_{1} C_2-C_1 S_2) \sqrt{{\cal J}}}{C^2_{1} + S^2_{1}}
\label{alphaR_1R_2}
\end{eqnarray}
where,
\begin{eqnarray}
{\cal J} = C^2_{1} + S^2_{1} - \chi^2_1.
\label{Lambda}
\end{eqnarray}
Here, ${\cal J}$ is analogous to ${\cal K}$ for the two centroid {\it ansatz}.
Eqs.(\ref{alphaBR_1}) and (\ref{alphaR_1R_2}) determine the centroid phase angles between the three clusters of the steady-state system. However, like 
Eq.(\ref{alphasol}) for two clusters, we can see that there are two conditions that will restore the time-dependence. 
Analogous to the role played by the sign of ${\cal K}$ in Eq.(\ref{alphasol}), the first condition is when ${\cal J}<0$ in 
Eqs.(\ref{alphaBR_1},\ref{alphaR_1R_2}) causing the right hand sides of 
Eqs.(\ref{alphaBR_1},\ref{alphaR_1R_2}) to become complex. The second condition, with no analog in the two cluster case, occurs 
when the modulus of the right hand side of Eq.(\ref{alphaR_1R_2}) is greater than unity, so that no solution exists for $\alpha_{R_1 R_2}$. 
Because of the many couplings, the threshold for this cannot be written in terms of a single
critical coupling.

We may determine when the fragmentation {\it ansatz} is preferable over that for the internally locked
scenario: if under a parameter variation, solutions to Eq.(\ref{alphaR_1R_2}) exist with $\alpha_{R_1 R_2}$
reasonably small, then {\it ansatz} Eq.(\ref{fixpoint1}) may be used over Eq.(\ref{fixpoint2}). Looking for the equivalent condition for $\alpha_{BR_1}$, means that the square of the
right hand side of  Eq.(\ref{alphaBR_1}) must be greater than one which, after some algebra, reduces to:
\begin{eqnarray*}
\left( \chi_1^2 - C_1^2 \right)^2 \left( C^2_1 +S^2_1 \right)<0,
\end{eqnarray*}
which may never hold.
A similar expression also holds for the internal locking scenario in Eq.(\ref{alphasol}). Hence a constant solution for $\alpha_{BR_1}$ will exist only if ${\cal J}>0$. Even if such a solution exists, there may be no solution for $\alpha_{R_1 R_2}$.

Bounds on stability are now straightforwardly determined, analogous to Eq.(\ref{stab-2clust}) using Taylor's `lemma' by contracting with the
zero eigenvectors of ${\cal L}_0'$, leading to the conditions:
\begin{eqnarray}
\cos(\phi-\alpha_{BR_1}) \geq  0, \;\;  \cos(\alpha_{R_1R_2}) \geq 0 , \nonumber \\
 \cos(\psi + \alpha_{BR_2}) \geq  - {{\sigma_R d^{(R_1 R_2)}_T}\over {\zeta_{RB} d^{(R_1 B)}_T}} \cos(\alpha_{R_1 R_2}).
\label{stab-3clust}
\end{eqnarray}
These three conditions appear to be mutually exclusive until we eliminate one of the phases using Eq.(\ref{13=12+23}), giving the
third condition:
\begin{eqnarray*}
 \cos(\psi + \alpha_{BR_1} + \alpha_{R_1 R_2}) &\geq & - {{\sigma_R d^{(R_1 R_2)}_T}\over {\zeta_{RB} d^{(R_1 B)}_T}} \cos(\alpha_{R_1 R_2}).
\end{eqnarray*}
As for the internal locking scenario, these conditions give a tighter (though still not definitive) bound on conditions for stability of the system than would Gershgorin's theorem. 
However, the symmetry of the two-dimensional $R_1 R_2$ block of the super-Laplacian Eq.(\ref{genLapl'})
suggests that the second condition of Eqs.(\ref{stab-3clust}), that on $\alpha_{R_1 R_2}$, is tight. Correspondingly,
when cross-interactions between Blue and $R_1$ are symmetric, and $\cos(\phi-\alpha_{BR_1})=\cos(\psi+\alpha_{BR_1})$
then all of ${\cal L}'$ is symmetric and all three conditions become tight bounds on stability.

\section{Measures of synchronisation}
To measure, on the one hand success in self-synchronisation within a given population, and on
the other hand success in one population achieving a certain phase shift with respect to the other, we use local order parameters
\begin{eqnarray*}
O_B = \frac{1}{N}\left|  \sum_{j\in {\cal B}} e^{i\beta_j} \right|, \;\; O_R = \frac{1}{M} \left| \sum_{j\in {\cal R}} e^{i\rho_j} \right|, O_{R_k} = \frac{1}{M_k} \left| \sum_{j\in {\cal R}_k} e^{i\rho_j} \right|,\;\; k \in \{1,2\}.
\end{eqnarray*}
The first two represent, for Blue and Red populations respectively, the order parameter introduced by \cite{Kura1984} for self-synchronisation: the closer $O_{B/R}$
is to the value one the more that the oscillators achieve phase synchronisation
{\it within their clusters}. 
The third, $O_{R_k}$, examines the phase synchronisation within the fragmented Red clusters in the case of three overall clusters. 
We emphasise that in all of our scenarios in the following the total system order parameter $1/(N+M) | \sum_{i\in{\cal B}} e^{i\beta_i} +\sum_{j\in{\cal R}} e^{i\rho_j}|$
will be far from the value one.

\section{Solution example: Hierarchy vs Random network }
Consider two populations of the same total number of agents but where a sub-population of each are cross-interacting
with same coupling and one-to-one:
\begin{eqnarray}
N = M,\;\; \zeta_{BR} = \zeta_{RB}, \label{simplify},
\end{eqnarray}
but with different internal networks and couplings and frustration parameters. 
Let the blue population form a hierarchy, naturally described by tree graphs, with a single root, and sets of distinct
branches connected to the root, each culminating in leaf nodes. Specifically we use a complete 4-kary tree, thus setting $N=21$.
We use a random Erd\"os-R\'enyi network 
for the Red population of the same number of agents generated by a link 
placed between nodes with probability $0.4$.
\begin{figure}
\begin{center}
\includegraphics[height=7cm]{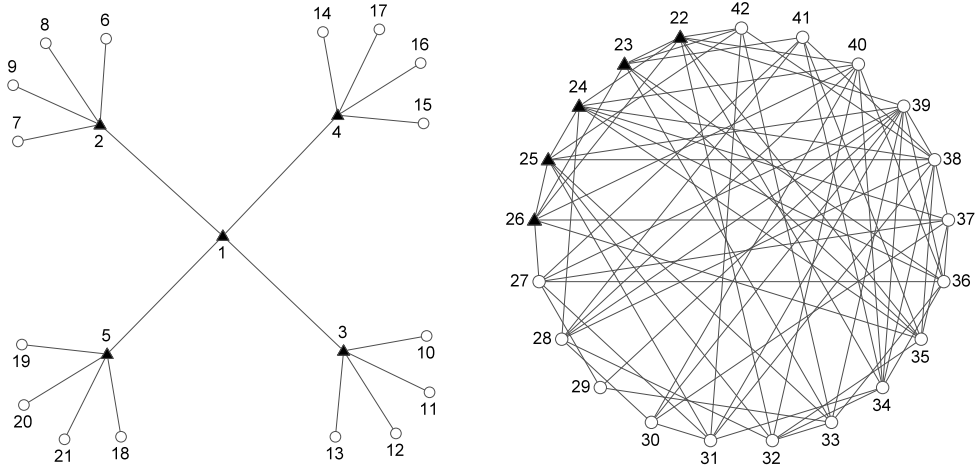}
\caption{Two networks on 21 nodes each, a hierarchical Blue (left) and random Red (right). Node $i$ in the tree graph with an open circle
interacts with a corresponding number, node $i+21$; all other nodes are represented by a filled triangle. } \label{twonets}
\end{center}
\end{figure}
The two networks are represented in Fig.\ref{twonets}.
We arrange the interaction between Blue and Red such that each leaf-node of Blue ($i=6,\dots,21$) interacts with the correspondingly labelled Red node ($i=27,\dots,42$), shown as open circles in Fig.\ref{twonets}; of course the random structure of Red means there is no consistent pattern in how those Red nodes link internally to other Red nodes 
(there is no hierarchically identifiable `leadership'). In other words,
$d_i^{(BR)}=d_i^{(RB)}= 1 \ \rm{or} \ 0$, for agents engaged, respectively not engaged, with a competitor, but $d^{(BR)}_T=16$.

\begin{figure}
\begin{center}
\includegraphics[height=4.5cm]{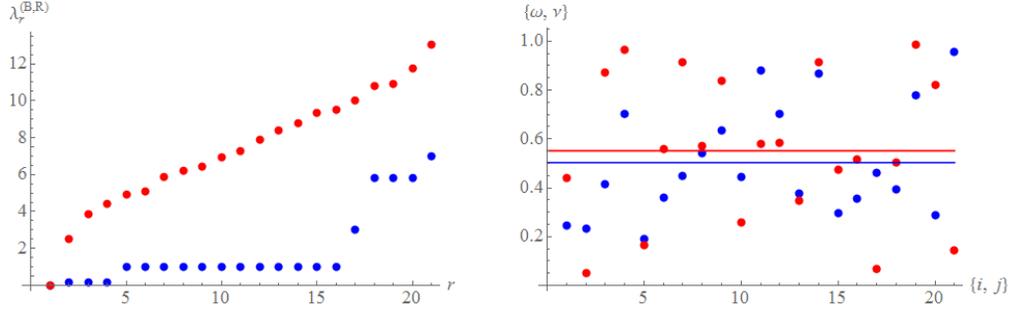}
\caption{Left: The spectrum of the graph Laplacian for Blue and Red networks coloured respectively blue and red. 
Right: The frequencies for Blue and Red agents according to the node labelling again coloured blue and red respectively, with solid lines indicating their corresponding means.} 
\label{spectra}
\end{center}
\end{figure}

The spectra of the respective graph Laplacians, and the frequencies of each agent, are shown in Fig.\ref{spectra}. A key observation to be made about the former 
(left, Fig.\ref{spectra}) are the many more low lying eigenvalues for
the Blue agents' network compared to that of Red - a 
consequence of the relative poor connectivity of a tree compared to a random graph.
The frequencies of each agent are drawn from a uniform distribution between
zero and one, $\omega_i,\nu_j\in [0,1]$. In the examples 
used here (right, Fig.\ref{spectra}), the average frequencies turn out to be $\bar{\omega}=0.503, \bar{\nu}=0.551$, so that if 
cross-couplings were set to zero the Red population would lap Blue over time. 

\subsection{Local phase synchronisation}
We choose the couplings as follows:
\begin{eqnarray}
\sigma_B=8, \;\; \sigma_R=0.5, \;\; \zeta_{BR}=\zeta_{RB}=0.4. \label{coupls}
\end{eqnarray}
As will be seen, these are sufficient to allow Blue and Red to achieve internal self-synchronisation in the absence of
cross-coupling. Note that $\sigma_B>\sigma_R$ is consistent both with tree networks being harder to synchronise given their relatively poor 
linking, and with the different profiles of low-lying Laplacian ($r=1,2,3)$ eigenvalues in Fig.\ref{spectra}.
In other words, the couplings are chosen such that $1/\sigma \lambda_r < 1$ as discussed for Eq.(\ref{critcond}).

\begin{figure}[t]
\begin{center}
\includegraphics[height=5.7cm]{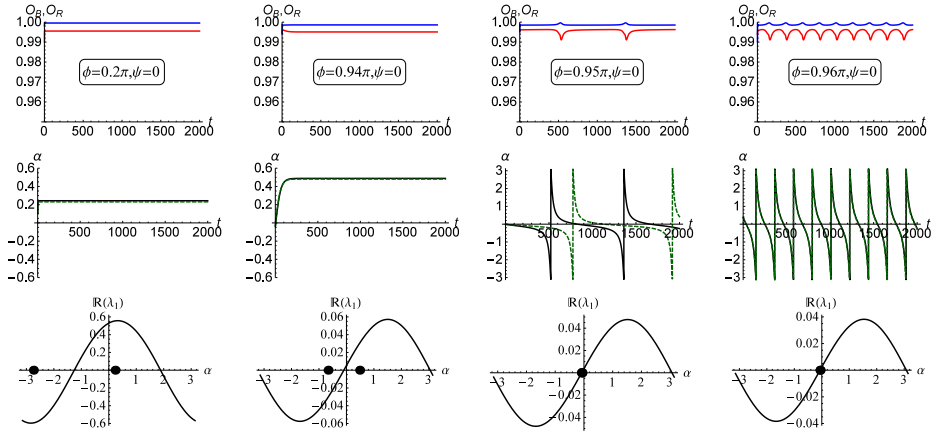}
\caption{
Behaviour of various measures with numerical solution of the full system Eqs.(\ref{Blue-eq},\ref{Red-eq}) 
with couplings in Eqs.(\ref{coupls}), $\psi=0$, and $\phi=0.2\pi$ (far left hand column), 
$\phi=0.94\pi$ 
(middle left column), $\phi= 0.95\pi$ (middle right hand column) and $\phi=0.96\pi$ (far right hand column).
Top row: plots of the local synchronisation order parameters $O_B$ (blue curve) and $O_R$ (red curve)
as functions of time for
the four values of $\phi$. 
Middle row: the average phase difference between Blue and Red $\alpha$ as a function of time for the four cases.
Also indicated is the time-dependent approximation for $\alpha$ of Eq.(\ref{t-alpha}) in the dashed line.
Bottom row: the real part of the lowest lying ($r=1$) eigenvalue of the super-Laplacian ${\cal L}$
in Eq.(\ref{genLapl}) as a function of $\alpha$ 
(all higher eigenvalues are strictly positive for all $\alpha$); the two solutions closest to the origin for $\alpha$ to Eq.(\ref{alpha-const-sol})
are shown as dots.
Note that only for $\phi=0.95\pi, 0.96\pi$ the solutions are complex and conjugates of each other so for these values there is an imaginary part
(identical for the two solutions). }
\label{lambdastudy}
\end{center}
\end{figure}

We examine numerical solutions to the full equations of motion  Eqs.(\ref{Blue-eq},\ref{Red-eq}) using 
Mathematica's NDSolve capability. We
fix $\psi=0$ but vary $\phi$ from zero upwards.
In Fig.\ref{lambdastudy} we plot a number of properties of the system for $\phi=0.2\pi, 0.94\pi, 0.95\pi, 0.96\pi$
as well as the behaviour of the lowest eigenvalue of ${\cal L}$ as a function of $\alpha$.

Running across the top row of plots in Fig.\ref{lambdastudy} we observe that
as $\phi$ crosses a threshold between $0.94\pi$ and $0.95\pi$ the Blue and Red systems
respectively change from near perfect phase synchrony $O_{B,R}\approx 1$ to a form of
cyclic synchrony. However the deviations from the value one are slight so that local phase locking is never significantly
destroyed; this was the key requirement for the linearisation approximation Eq.(\ref{approx}) to hold.
The middle row of figures provide further insight. The average phase difference $\alpha$
agrees with the approximation in Eq.(\ref{t-alpha}) and is constant, up to
$\phi=0.94\pi$, confirming that Blue and Red are globally phase locked, with Blue only slightly ahead of
Red. However at $\phi=0.95\pi$
the apparent locking is very temporary - the Blue population initially slips behind Red (seen in $\alpha$ 
becoming negative), rotates through $2\pi$ until once again the system temporarily locks.
The same change in behaviour can be seen in the approximate solution for $\alpha(t)$ 
due to ${\cal K}$ becoming negative in the $\tanh$ of Eq.(\ref{t-alpha}) though the period disagrees at $\phi=0.95\pi$,
a consequence of being very close to the critical point. As $\phi$ increases further, the numerical and approximate analytical results for $\alpha$ agree more
precisely.

\begin{figure}[t]
\begin{center}
\includegraphics[height=5cm]{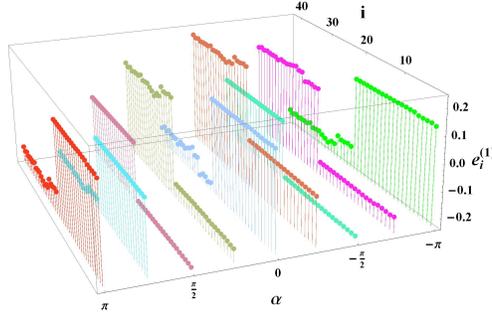}
\caption{Value of the components $i$ of the eigenvector corresponding to the lowest non-zero eigenvalue of
${\cal L}$ as a function of $\alpha$.} \label{evec}
\end{center}
\end{figure}

We turn next to the behaviour of the lowest eigenvalue of ${\cal L}$ and the solutions
to Eq.(\ref{alpha-const-sol}). For $\phi\leq 0.94\pi$ there are two lowest solutions to Eq.(\ref{alpha-const-sol}), both real. The bottom left and middle
plots of Fig.\ref{lambdastudy} show the dependence of the lowest eigenvalue of ${\cal L}$ on
$\alpha$, $\lambda_1(\alpha)$,  as expressed through Eq.(\ref{genLapl}). We state (but do not show in the plot) that for $r>1$
the eigenvalues turn out to be positive (and real) for all $\alpha$. That the single eigenvalue for $r=1$ has 
the potential to change sign is associated with the structure of the eigenvector. 
We see in  Fig.\ref{evec} that $e^{(1)}_i$ has a step like 
dependence at $i=N,N+1=21,22$ where the index for Blue nodes changes to those for Red.
This is seen across various values of $\alpha$.
Thus, the associated eigenvector distinguishes only the Blue-Red partition of the overall network, and not finer structures within.
In other words, the instability will only emerge for the dynamics of Blue in relation to Red and not for any substructures within them.

Returning to the bottom far left and middle left plots of Fig.\ref{lambdastudy}, 
for the two lowest values of $\alpha$ which solve Eq.(\ref{alpha-const-sol}),
the lowest Laplacian eigenvalue is respectively negative (left-most
dot) and positive (right-most dot); thus for the right-most solution for $\alpha$
{\it the lowest normal mode is
stable}. 
The crossing of the eigenvalue at zero is close to, but not exactly, the
condition $\alpha=\phi\pm \frac{\pi}{2}$ from the first of the conditions Eq.(\ref{stab-2clust}).

At $\phi=0.95\pi$ the solution to  Eq.(\ref{alpha-const-sol})
is complex: the two solutions are conjugates of each other so that the real parts are identical.
Thus in the bottom middle and right plots of the eigenvalue only one point appears - the common real part to the two solutions.
This point corresponds to ${\cal K}\rightarrow 0^{-}$ so it gives the threshold for
time-dependence in  Eq.(\ref{specialK}). On the other hand, with $\psi=0$ and the Blue-Red interactions symmetric in this example
then the super-Laplacian is approximately symmetric. The two conditions Eq.(\ref{stab-2clust}) approximately give the same 
condition for a zero eigenvalue which in turn agree with the crossing at zero of the curve in Fig.\ref{lambdastudy}.
Thus at this point the zero mode and normal mode dynamics are linked: the onset of dynamics in the leading order solution to $\alpha$
coincides with the appearance of a negative eigenvalue of the super-Laplacian.

Concluding this section we find there is a regime of the full system dynamics in which our linearisation 
approximation with locking internal to the populations holds showing behaviours consistent with the analytic solutions. The fixed points may be 
approximated through Eq.(\ref{alpha-const-sol}), where the stable one can be identified by
the behaviour of the spectrum of ${\cal L}$. In the context of the Blue population seeking a significantly large frustration value,
to operate half a cycle ahead of the competitor ($\phi\rightarrow\pi$), too much is being sought.
The dynamics
only permit a small advantage to Blue in the small positive solutions $\alpha$ 
consistent with stability, despite the increasing value of $\phi$. However, beyond the critical point
the dynamics yields the advantage to Red - albeit temporarily.
The Blue population is indeed beyond the point of diminishing return on its strategy.

\subsection{Optimal strategy for Blue for internally locked scenario}
\begin{figure}[t]
\begin{center}
\includegraphics[height=5cm]{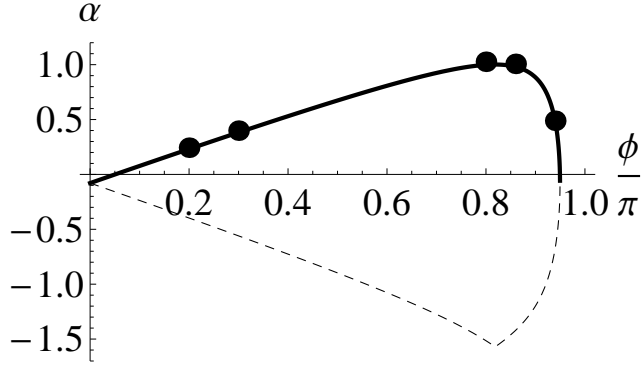}
\caption{
Solution $\alpha(\phi,\psi)$ from the $\arcsin$ of Eq.(\ref{alphasol}) as a function of $\phi/\pi$
with $\psi=0$ and couplings given by Eq.(\ref{coupls}): 
solid line is the negative choice of the root, dashed line is the positive root, and the points are the values of 
the numerical $\alpha$ at $t=2000$ at $\phi/\pi=0.2, 0.3, 0.8, 0.86, 0.94$.
 }
\label{optimalph}
\end{center}
\end{figure}

This point of diminishing return for Blue may be analytically computed through
Eqs.(\ref{alpha-const-sol},\ref{alphasol}) as part of the optimisation strategy discussed at that stage.
We plot in Fig.\ref{optimalph} the solutions for $\alpha(\phi,\psi)$ from the $\arcsin$ of Eq.(\ref{alphasol}): the
solid line is the negative root and the dashed line the positive root of the quadratic equation Eq.(\ref{alphasol}). The corresponding plot for $\sin \alpha$ only differs around the turning points. We superimpose the numerical steady state values
of $\alpha(t)$ at $t=2000$ to check which of the roots of Eq.(\ref{alphasol}) is the relevant solution.

We observe a steady rise in the value of $\alpha$ and then a turning point and asymptote down: this is the point at
which an imaginary part develops and there is no longer a constant solution to Eq.(\ref{alphadot-approx}).
The turning point occurs at $\phi=0.82\pi$, at which point $\sin\alpha=0.84$ or $\alpha=1.00=0.32\pi$.
Thus the optimal strategy for the Blue population is for its agents interacting with competitors to seek to be $\phi=0.82\pi$,
though the effects of interactions will limit the {\it actual lead} to approximately
one sixth of a cycle. 

\subsection{Fragmentation regime}
\begin{figure}[t]
\begin{center}
\includegraphics[height=5.5cm]{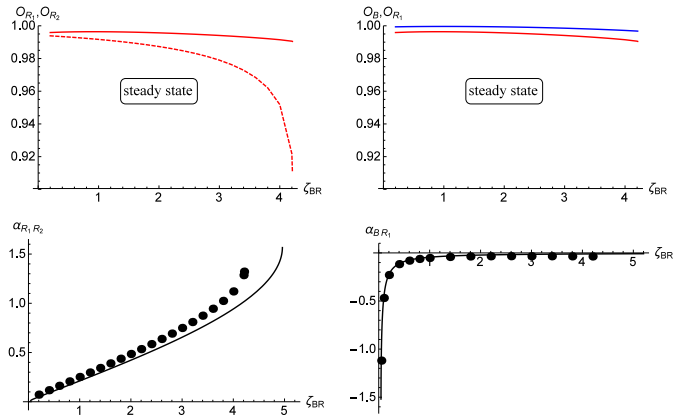}
\caption{
Behaviour of various measures with numerical solution of the full system Eqs.(\ref{Blue-eq},\ref{Red-eq}) 
with couplings in Eqs.(\ref{coupls2}) but varying $\zeta_{BR}=\zeta_{RB}$.
Top row: plots of the local synchronisation order parameters $O$ for the two Red populations,
${\cal R}_1$ which interacts with Blue (left panel, solid red), and ${\cal R}_2$ (left panel, dashed red curve);
and for the Blue population (right panel, blue curve) and the Red population ${\cal R}_1$ with which it interacts (right panel, solid red curve).
Bottom row: the angle between the centroids within Red (left) and between Blue and ${\cal R}_1$ (right), with dots the numerical
solution and the solid curve the analytical approximation 
 Eqs.(\ref{alphaR_1R_2},\ref{alphaBR_1}).
}
\label{KBRstudy}
\end{center}
\end{figure}

We now show behaviours of the full system but where our requirement for local phase synchronisation
ceases to hold and fragmentation within a population occurs.
This time, we select $\zeta_{BR}=\zeta_{RB}$ and parameters:
\begin{eqnarray}
\sigma_B=8, \;\; \sigma_R=0.5, \;\; \psi=\phi=\pi/4. \label{coupls2}
\end{eqnarray}
Thus we have Blue and Red both seeking a significant (quarter of a cycle) and identical advantage over the other.
In Fig.\ref{KBRstudy} we plot local order parameters at large time $(t=2000)$ for increasing $\zeta_{BR}$.
We also plot the angle between centroids of various sub-populations, comparing both numerical solutions and 
analytical solutions.

The top row of plots in Fig.\ref{KBRstudy} show that as $\zeta_{BR}$ increases 
the part of Red interacting with Blue, ${\cal R}_1$, maintains good synchronisation (solid red curve) while the other part, ${\cal R}_2$, suffers increasingly lower levels of synchronisation
(dashed red curve); it may be described as approaching a `splay state' \cite{Nord2015} with individual oscillator phases spread across the unit circle. The Blue population however also maintains good synchronisation (solid blue curve, right panel).
Thus Blue and the parts of Red with which its agents interact enjoy increased locking with respect to each other at the cost of Red's internal coherence;
Red not only fragments into two parts, but those agents not interacting with Blue become increasingly incoherent. 
In the bottom two panels of Fig.\ref{KBRstudy} we compare the angle between the clusters as numerically determined, the solid dots, and our analytical 
solutions. We see firstly that the angle between the two parts of Red increases with $\zeta_{BR}$, consistent with their growing separation.
However our approximation, the solid curve, increasingly deviates from the numerical result for $\alpha_{R_1R_2}$ - the former
is predicated on the fragmented sub-population retaining some level of internal synchronisation.
On the other hand, in the right hand panel we see solid agreement between analytical and numerical solutions for $\alpha_{BR_1}$ 
across values of $\zeta_{BR}$, with the angle between Blue agents and their competitors diminishing with
increasing interaction strength.

\begin{figure}
\begin{center}
\includegraphics[height=3cm]{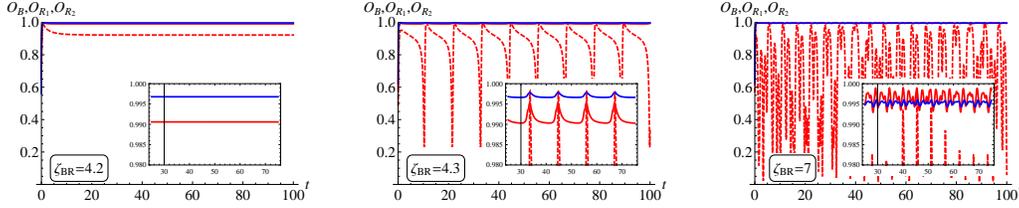}
\caption{
Time dependence of the order parameters for Blue, $O_B$, (solid blue curve) and the two Red fragments, $O_{R_1}, O_{R_2}$, (solid and dashed red curves, respectively) 
at various cross-couplings 
$\zeta_{BR}$ at and beyond the critical point; insets show the behaviour in the vicinity of value one.
}
\label{evaporation}
\end{center}
\end{figure}

In Fig.\ref{KBRstudy} we plot only up to the point where $\alpha_{R_1R_2}$ remains time-independent. For $\zeta_{BR}>4.2$,
$\alpha_{R_1R_2}$ becomes time-dependent, notably when it is close to $\pi/2$, where Eq.(\ref{stab-3clust}) predicts 
negative eigenvalues of the super-Laplacian ${\cal L}'$. This time-dependence is shown in Fig.\ref{evaporation} where the first plot gives the behaviour
of the order parameters at $\zeta_{BR}=4.2$ where the points in Fig.\ref{KBRstudy} finish, and the subsequent two plots show higher values of $\zeta_{BR}$.
For example, at $\zeta_{BR}=4.3$ (middle plot, Fig.\ref{evaporation}) strong cyclic behaviour manifests in $O_{R_2}$ showing that the second fragment of Red has itself fragmented in
two. The inset in this case zooms in on the time-dependence near value one showing that Blue and the Red fragment locked with it show very small fluctuations as a 
consequence of the interactions wth Red. At higher $\zeta_{BR}$ the behaviour becomes more complicated though still periodic with second and third order subcycles in the 
${\cal R}_2$ fragment consistent with multiple subfragments. Beyond this value of $\zeta_{BR}$ no coherence remains in ${\cal R}_2$ and it may be said to have `evaporated' or `disintegrated'.
We emphasise that though we are unable to describe analytically this chaotic region, we are able to determine the threshold before such behaviour breaks out - our
key objective in this paper.

\begin{figure}
\begin{center}
\includegraphics[height=3.5cm]{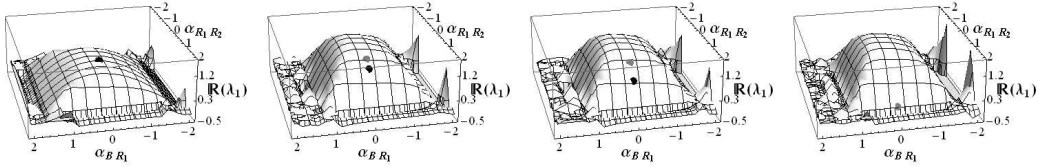}
\caption{
Plot of the lowest eigenvalue of the super-Laplacian as a function of $\alpha_{BR_1}$ and $\alpha_{R_1R_2}$
for four values of $\zeta_{BR}=1,4,4.2,5$. The mesh is chosen to be finer for values less than or equal to zero. Superimposed 
as dots are the numerical (black) and analytical (grey) solutions for the angles. For the right-most plot, there is only a grey dot as numerically the
angles $\alpha_{BR_1}$ and $\alpha_{R_1R_2}$ are time-dependent, while Eqs.(\ref{alphaBR_1},\ref{alphaR_1R_2})
still give a steady-state solution.
}
\label{KBRevals}
\end{center}
\end{figure}

We plot next the $r=1$ eigenvalue of the super-Laplacian as a function of the two angles $\alpha_{BR_1}$ and $\alpha_{R_1R_2}$
in Fig.\ref{KBRevals} for a range of values of $\zeta_{BR}$.
We superimpose on these plots the positions of the numerically and analytically determined 
(from Eqs.(\ref{alphaBR_1},\ref{alphaR_1R_2})) centroid angles for the various values of $\zeta_{BR}$.
Note that we do not plot for the full range $(-\pi,\pi)$ - we have verified that over the range given in  Fig.\ref{KBRevals},
the $r=1$ eigenvalue is the {\it lowest} and is the {\it first} to cross zero; outside of the range shown eigenvalues cross and the shape
becomes considerably complicated.

We observe firstly that the eigenvalues are consistently positive close to the bounds determined from the first two of
Eqs.(\ref{stab-3clust}); indeed by inspection from Fig.\ref{KBRevals}, the eigenvalue crosses zero for $\alpha_{R_1 R_2}=\pm\frac{\pi}{2}$,
exactly in accordance with the second of Eqs.(\ref{stab-3clust}), consistent with that condition being a tight bound on the spectrum.
We see that up to $\zeta_{BR}=4.2$ the numerical solution for the angles $\alpha_{BR_1}$ and $\alpha_{R_1R_2}$ sits
on the surface, consistent with there being a positive eigenvalue and thus a stable fixed point. However, with increasing $\zeta_{BR}$
there is a divergence between the numerical behaviour and analytically determined solutions in the $\alpha_{R_1 R_2}$
directions - consistent with Fig.\ref{KBRstudy}, and the fact that synchronisation is breaking down within the Red population.
Finally, at $\zeta_{BR}=5$ there is no longer steady-state behaviour for $\alpha_{BR_1}$ and $\alpha_{R_1R_2}$, while
Eqs.(\ref{stab-3clust}) still give such a solution. Thus only a grey dot is indicated in the right-most panel of  Fig.\ref{KBRevals}.
However this value exactly coincides with a vanishing of the $r=1$ eigenvalue of ${\cal L}'$.
As for the internally locked scenario, the point at which a negative eigenvalue appears coincides with the point at which
the angles between centroids of clusters become time-dependent, with Eq.(\ref{alphaR_1R_2}) no longer having a solution.
Indeed, as $\zeta_{BR}$ is increased with all else fixed, ${\cal J}$ increases, as seen from Eq.(\ref{Lambda}), so 
that variations in this way never cross ${\cal J}=0$, consistent with the fragmentation {\it ansatz} being the appropriate
fixed point choice.

\subsubsection{Optimal strategy for Blue under Red fragmentation}
An `optimal strategy' for Blue in this case is more context dependent. If it is more important that the competitor remain coherent
and Blue maintain a lead in the cycle, then the bottom two plots of Fig.\ref{KBRstudy} offer a guide. With 
the conflicting strategies of $\phi=\psi=\frac{\pi}{4}$, the best that Blue can hope for is to dynamically match Red
with $\alpha_{BR_1}=0$ at large $\zeta_{BR}$. However, beyond $\zeta_{BR}=4.2$ there is diminishing return.
This point can still be reasonably estimated analytically up to about $\zeta_{BR}=4$. However, if complete dislocation of the competitor
is desired then $\zeta_{BR}>4.2 $ is desirable - and this bound too can be estimated from the 
point where the analytical result Eq.(\ref{alphaR_1R_2}) fails to have a solution.
Indeed, since the behaviour in this regime exhibits Red agents incoherent with respect to each other {\it and} Blue agents tightly locked
with their Red competitor agents, this is truly a case of Blue having its cake and eating it!

\section{Conclusions}

We have formulated a mathematical model, based on the two-network frustrated Kuramoto model, for two populations, labelled Blue and Red, on separate networks interacting 
internally in a cooperative fashion, to self-synchronise individual limit cycles via zero frustration parameters, and externally
in a competitive manner with non-zero frustration parameters, with each seeking to be some phase difference in relation to the limit cycle of the other. 
We have determined threshold conditions for two types of
fixed point behaviour of this combined system, one with both sides achieving their respective internal goals and one or both sides achieve their external goal,
and another with one side reaching a threshold for failure to achieve the internal goal.
In both cases,
we could reduce the criterion for a stable fixed point to a closed form analytical result from which optimal choices of 
frustration parameters may be computed.
The stability of the fixed point relies on positivity of eigenvalues of a generalisation of the graph
Laplacian. We examined some specific cases, comparing full numerical solutions to the analytic solutions, 
both where our local phase synchronisation assumption is respected for increasing couplings, and where it fails. In 
both cases, the analytic formulae successfully predict the behaviour of the 
full non-linear system, including the critical point at which Blue and Red fail to lock with respect to each other
or one internally fragments. In the latter cases, the agreement is poorer because of the increasing
`splaying' and `dislocation' in the fragmenting population prior to time-dependence kicking in. Again, in both cases we demonstrated
how optimisation of the frustration may be computed by determination of thresholds which guarantee avoidance of `uncontrollable' chaotic behaviour of parts of the
system.

We have shown an example of optimisation in this paper only to illustrate the value of our mathematical approach but have 
deliberately avoided searching at this stage for Nash equilibria. Future work will address the impact of noise and the extension to
multiple internal fragmentation and disintegration of clusters. Our hope is in the long term to 
provide mathematical tools to guide the design of networks to structure against known - or at least statistically parametrisable - competitors.

% Appendix here
% Options are (1) APPENDIX (with or without general title) or
%             (2) APPENDICES (if it has more than one unrelated sections)
% Outcomment the appropriate case if necessary
%
% \begin{APPENDIX}{<Title of the Appendix>}
% \end{APPENDIX}
%
%   or
%
% \begin{APPENDICES}
% \section{<Title of Section A>}
% \section{<Title of Section B>}
% etc
% \end{APPENDICES}

%%
%\theendnotes

% Acknowledgments here
\section*{Acknowledgements}
The authors gratefully acknowledge discussions with 
Richard Taylor, Tony Dekker, Iain Macleod and Markus Brede.
ACK is supported through a Chief Defence Scientist Fellowship.

\end{document}